\documentclass[12pt]{article}%
\usepackage{graphicx}
\usepackage[intlimits]{amsmath}
\usepackage{latexsym}
\usepackage{amsfonts}
\usepackage{amssymb}%
\makeatletter
\newfont{\footsc}{cmcsc10 at 8truept}
\newfont{\footbf}{cmbx10 at 8truept}
\newfont{\footrm}{cmr10 at 10truept}
\newtheorem{theorem}{Theorem}

\newtheorem{proposition}[theorem]{Proposition}

\newenvironment{proof}[1][Proof]{\noindent{\textbf {#1}  }}  {\hfill$\Box$\bigskip}

\begin{document}

\title{A contribution to the Zarankiewicz problem}
\author{Vladimir Nikiforov\\Department of Mathematical Sciences, University of Memphis, \\Memphis, TN 38152, USA, email: \textit{vnikifrv@memphis.edu}}
\maketitle

\begin{abstract}
Given positive integers $m,n,s,t,$ let $z\left(  m,n,s,t\right)  $ be the
maximum number of ones in a $\left(  0,1\right)  $ matrix of size $m\times n$
that does not contain an all ones submatrix of size $s\times t$. We show that
if $s\geq2$ and $t\geq2,$ then for every $k=0,\ldots,s-2,$%
\[
z\left(  m,n,s,t\right)  \leq\left(  s-k-1\right)  ^{1/t}nm^{1-1/t}+kn+\left(
t-1\right)  m^{1+k/t}.
\]
This generic bound implies the known bounds of K\"{o}vari, S\'{o}s and
Tur\'{a}n, and of F\"{u}redi. As a consequence, we also obtain the following results:

Let $G$ be a graph of $n$ vertices and $e\left(  G\right)  $ edges, and let
$\mu$ be the spectral radius of its adjacency matrix. If $G$ does not contain
a complete bipartite subgraph $K_{s,t},$ then the following bounds hold%
\[
\mu\leq\left(  s-t+1\right)  ^{1/t}n^{1-1/t}+\left(  t-1\right)
n^{1-2/t}+t-2,
\]
and
\[
e\left(  G\right)  <\frac{1}{2}\left(  s-t+1\right)  ^{1/t}n^{2-1/t}+\frac
{1}{2}\left(  t-1\right)  n^{2-2/t}+\frac{1}{2}\left(  t-2\right)  n.
\]

\textbf{Keywords: }\textit{bipartite subgraphs; Zarankiewicz problem; spectral
radius. \smallskip}


\end{abstract}

\subsection*{Introduction}

How large can be the spectral radius $\mu$ of a graph order $n$ that does not
contain a complete bipartite subgraph $K_{s,t}?$ This is a spectral version of
the famous Zarankiewicz problem: \emph{how many edges can have a graph of
order }$n$\emph{ if it does not contain }$K_{s,t}?$ Except for few cases, no
satisfactory solution to either of these problems is known. In an unpublished
pioneering work, Babai and Guiduli (see, e.g., \cite{Gui96}) have shown that
\[
\mu\leq\left(  \left(  s-1\right)  ^{1/t}+o\left(  1\right)  \right)
n^{1-1/t}.
\]
Using a different method, here we improve this result as follows:

\begin{theorem}
\label{th1}Let $s\geq t\geq2,$ and let $G$ be a $K_{s,t}$-free graph of order
$n$ and spectral radius $\mu.$ If $t=2,$ then%
\begin{equation}
\mu\leq1/2+\sqrt{\left(  s-1\right)  \left(  n-1\right)  +1/4}.\label{in0}%
\end{equation}
If $t\geq3,$ then%
\begin{equation}
\mu\leq\left(  s-t+1\right)  ^{1/t}n^{1-1/t}+\left(  t-1\right)
n^{1-2/t}+t-2.\label{in1}%
\end{equation}

\end{theorem}

Below we show that the bounds (\ref{in0}) and (\ref{in1}) are tight for some
values of $s$ and $t.$ On the other hand, in view of the inequality $2e\left(
G\right)  \leq\mu n,$ we see that if $G$ is a $K_{s,t}$-free graph of order
$n,$ then
\begin{equation}
e\left(  G\right)  \leq\frac{1}{2}\left(  s-t+1\right)  ^{1/t}n^{2-1/t}%
+\frac{1}{2}\left(  t-1\right)  n^{2-2/t}+\frac{1}{2}\left(  t-2\right)
n.\label{in2}%
\end{equation}
This is a slight improvement of a result of F\"{u}redi \cite{Fur96}.

To prove Theorem \ref{th1}, we first find a family of new upper bounds for the
matrix Zarankiewicz problem, thereby extending some previous results.

\subsection*{The matrix Zarankiewicz problem}

Let $J_{s,t}$ denote the all ones matrix of size $s\times t.$ Given positive
integers $m,n,s,t,$ let $z\left(  m,n,s,t\right)  $ be the maximum number of
ones in a $\left(  0,1\right)  $ matrix of size $m\times n$ that does not
contain $J_{s,t}$ as a submatrix. 

Here is an equivalent definition: $z\left(  m,n,s,t\right)  $ is the maximum
number of edges in a bipartite graph $G$ with vertex classes $A$ of size $n$
and $B$ of size $m$ such that $G$ does not contain a copy of $K_{s,t}$ with
vertex class of size $s$ in $A$ and vertex class of size $t$ in $B.$

The problem of finding $z\left(  m,n,s,t\right)  $ is known as the general
Zarankiewicz problem. In \cite{KST54}, K\"{o}vari, S\'{o}s and Tur\'{a}n gave
one of the earliest bounds on $z\left(  m,n,s,t\right)  ,$ which in simplified
form reads as
\begin{equation}
z\left(  m,n,s,t\right)  \leq\left(  s-1\right)  ^{1/t}nm^{1-1/t}+\left(
t-1\right)  m.\label{KSTin}%
\end{equation}
Later, F\"{u}redi \cite{Fur96} improved this bound showing that if $s\geq t,$
then
\begin{equation}
z\left(  m,n,s,t\right)  \leq\left(  s-t+1\right)  ^{1/t}nm^{1-1/t}%
+tm^{2-2/t}+tn.\label{Furin}%
\end{equation}
The proof of F\"{u}redi, although rather involved, is based on double counting
as in \cite{KST54}. Using a different approach, we show that, in fact,
(\ref{Furin}) and (\ref{KSTin}) are particular cases of a whole sequence of
subtler bounds on $z\left(  m,n,s,t\right)  $. Instead of using double
counting, we start with (\ref{KSTin}) and deduce by induction a number of
inequalities, one of which implies (\ref{Furin}). The following theorem gives
the precise statement.

\begin{theorem}
\label{thZ}If $s\geq2$ and $t\geq2,$ then for every $k=0,\ldots,s-2,$%
\begin{equation}
z\left(  m,n,s,t\right)  \leq\left(  s-k-1\right)  ^{1/t}nm^{1-1/t}+\left(
t-1\right)  m^{1+k/t}+kn.\label{ourin}%
\end{equation}

\end{theorem}

Given (\ref{ourin}), letting $k=0,$ we obtain the bound of K\"{o}vari, S\'{o}s
and Tur\'{a}n (\ref{KSTin}). Also, if $s\geq t,$ letting $k=t-2,$ we obtain
\[
z\left(  m,n,s,t\right)  <\left(  s-t+1\right)  ^{1/t}nm^{1-1/t}+\left(
t-1\right)  m^{2-2/t}+\left(  t-2\right)  n,
\]
which is a slight improvement of F\"{u}redi's bound (\ref{Furin}). 

At first glance it is unclear whether the parameter $k$ is really useful in
inequality (\ref{ourin}). Indeed, for $n=m,$ setting $k=\min\left\{
s,t\right\}  -2$ gives the best inequality for $n$ large enough. However, for
arbitrary $n$ and $m,$ the parameter $k$ can give additional improvement, as
shown in the following proposition, whose proof is omitted.

\begin{proposition}
Let $s\geq3$ and $t\geq3,$ and $0\leq k\leq s-2.$ There exist $A=A\left(
s,t,k\right)  >0$ and $B=B\left(  s,t,k\right)  >0$ such that for all
sufficiently large $n$ and $m$ satisfying
\[
Am^{\left(  k+1\right)  /t}\leq n\leq Bm^{\left(  k+2\right)  /t},
\]
we have%
\[
\left(  s-k-1\right)  ^{1/t}nm^{1-1/t}+\left(  t-1\right)  m^{1+k/t}%
+kn<\left(  s-i-1\right)  ^{1/t}nm^{1-1/t}+\left(  t-1\right)  m^{1+i/t}+in
\]
for all $i\in\left[  0,s-2\right]  \backslash\left\{  k\right\}  .$
\end{proposition}

\subsection*{Tightness of the bounds (\ref{in0}) and (\ref{in1})}

For some values of $s$ and $t$ the bounds given by (\ref{in0}) and (\ref{in1})
are tight.

\subsubsection*{The case $t=2$}

For $s=t=2$ inequality (\ref{in0}) gives that every $K_{2,2}$-free graph $G$
of order $n$ satisfies
\[
\mu\left(  G\right)  \leq1/2+\sqrt{n-3/4}.
\]
This bound is tight: equality holds for the friendship graph. Note that
letting $q$ be a prime power, the Erd\H{o}s-Renyi polarity graph is a
$K_{2,2}$-free graph of order $n=q^{2}+q+1$ and $q\left(  q+1\right)  ^{2}/2$
edges. Thus, its spectral radius $\mu\left(  ER_{q}\right)  $ satisfies
\[
\mu\left(  ER_{q}\right)  \geq\frac{q^{3}+2q^{2}+q}{q^{2}+q+1}>q+1-\frac{1}%
{q}=1/2+\sqrt{n-3/4}-\frac{1}{\sqrt{n}-1},
\]
which is also close to the upper bound.

For $s>2,$ equality  in (\ref{in0}) is attained when $G$ is a strongly regular
graph in which every two vertices have exactly $s-1$ common neighbors. There
are examples of strongly regular graphs of this type; here is a small
selection from Gordon Royle's webpage:%

\[%
\begin{tabular}
[c]{||c|c|c||}\hline\hline
\emph{s} & \emph{n} & $\mu\left(  G\right)  $\\
$3$ & $45$ & $12$\\
$4$ & $96$ & $20$\\
$5$ & $175$ & $30$\\
$6$ & $36$ & $15$\\\hline\hline
\end{tabular}
\]
We are not aware whether there are infinitely many strongly regular graphs in
which every two vertices have the same number of common neighbors. However,
F\"{u}redi \cite{Fur96a} has shown that for any $n$ there exist $K_{s,2}$-free
graph $G_{n}$ of order $n$ such that%
\[
e\left(  G_{n}\right)  \geq\frac{1}{2}n\sqrt{sn}+O\left(  n^{4/3}\right)  ,
\]
and so,%
\[
\mu\left(  G_{n}\right)  \geq\sqrt{sn}+O\left(  n^{1/3}\right)  ;
\]
thus (\ref{in0}) is tight up to low order terms.

\subsubsection*{The case $s=t=3$}

The bound (\ref{in1}) implies that if $G$ is a $K_{3,3}$-free graph of order
$n,$ then
\[
\mu\left(  G\right)  \leq n^{2/3}+2n^{1/3}+1\text{\textbf{.}}%
\]

On the other hand, a construction due to Alon, R\`{o}nyai and Szab\`{o}
\cite{ARS99} implies that for all $n=q^{3}-q^{2},$ where $q$ is a prime power,
there exists a $K_{3,3}$-free graph $G_{n}$ of order $n$ with
\[
\mu\left(  G_{n}\right)  \geq n^{2/3}+\frac{2}{3}n^{1/3}+C
\]
for some constant $C>0.$ Thus, the bound (\ref{in1}) is asymptotically tight
for $s=t=3.$ The same conclusion can be obtained from Brown's construction of
$K_{3,3}$-free graphs \cite{Bro66}.

\subsubsection*{The general case}

As proved in \cite{ARS99}, there exists $c>0$ such that for all $t\geq2$ and
$s\geq\left(  t-1\right)  !+1$, there is a $K_{s,t}$-free graph $G_{n}$ of
order $n$ with
\[
e\left(  G_{n}\right)  \geq\frac{1}{2}n^{2-1/t}+O\left(  n^{2-1/t-c}\right)  .
\]
Hence, for such $s$ and $t$ we have%
\[
\mu\left(  G\right)  \geq n^{1-1/t}+O\left(  n^{1-1/t-c}\right)  ;
\]
thus, the bound (\ref{in1}) and the earlier bound of Babai and Guiduli give
the correct order of the main term.

\subsection*{Proof of Theorem \ref{thZ}}

\paragraph{Some matrix notation}

Let $\left\vert X\right\vert $ denote the cardinality of a finite set $X.$ Let
$A=\left(  a_{ij}\right)  $ be a $\left(  0,1\right)  $-matrix, and let the
rows and columns of $A$ be indexed by the elements of two disjoint sets
$R\left(  A\right)  $ and $C\left(  A\right)  .$ Then:

- for any $i\in R,$ we let $C_{i}=\left\{  j:j\in C\left(  A\right)  ,\text{
}a_{ij}=1\right\}  $ and set $r_{i}=\left\vert C_{i}\right\vert ;$

- for any $j\in C,$ we let $R_{i}=\left\{  i:i\in R\left(  A\right)  ,\text{
}a_{ij}=1\right\}  $ and set $c_{j}=\left\vert R_{j}\right\vert ;$

- $\left\Vert A\right\Vert $ stands for the sum of the entries of $A;$

- given nonempty sets $I\subset R\left(  A\right)  ,$ $J\subset C\left(
A\right)  ,$ we write $A\left[  I,J\right]  $ for the submatrix of the entries
$a_{ij}$ satisfying $i\in I,$ $j\in J$.\bigskip

\begin{proof}
[\textbf{Proof of Theorem \ref{thZ}}]We shall use induction on $k.$ For $k=0$,
the assertion is given by (\ref{KSTin}). Suppose $k\geq1$ and assume the
assertion true for all $k^{\prime}<k.$ Let $A=\left(  a_{ij}\right)  $ be a
$\left(  0,1\right)  $-matrix of size $m\times n,$ and let $R=R\left(
A\right)  ,$ $C=C\left(  A\right)  .$ Suppose that $A$ does not contain
$J_{s,t}$ as a submatrix and that $k\leq s-2$. Our goal is to prove that%
\[
\left\Vert A\right\Vert \leq\left(  s-k-1\right)  ^{1/t}nm^{1-1/t}+\left(
t-1\right)  m^{1+k/t}+kn.
\]
Select $i\in R$ and define the sets
\[
U=R\backslash\left\{  i\right\}  ,\ \ \ W=C_{i}.
\]
Note that the matrix $A\left[  U,W\right]  $ does not contain $J_{s-1,t}$ as a
submatrix since the $i$'th row of $A\left[  R,W\right]  $ consists of all ones
and we would have a $J_{s,t}$ in $A.$ Therefore,
\[
\left\Vert A\left[  U,W\right]  \right\Vert \leq z\left(  \left\vert
U\right\vert ,\left\vert W\right\vert ,s-1,t\right)  ,
\]
and by the induction assumption applied for $s-1$ and $k-1,$ we have
\begin{align}
\left\Vert A\left[  U,W\right]  \right\Vert  &  \leq\left(  s-k-1\right)
^{1/t}\left\vert W\right\vert \left\vert U\right\vert ^{1-1/t}+\left(
t-1\right)  \left\vert U\right\vert ^{1+\left(  k-1\right)  /t}+\left(
k-1\right)  \left\vert W\right\vert \nonumber\\
&  \leq\left(  s-k-1\right)  ^{1/t}r_{i}m^{1-1/t}+\left(  t-1\right)
m^{1+\left(  k-1\right)  /t}+\left(  k-1\right)  r_{i}.\label{in3}%
\end{align}
A closer look at $A\left[  U,W\right]  $ shows that
\begin{align*}
\left\Vert A\left[  U,W\right]  \right\Vert  &  =\sum_{j\in C_{i}}\sum_{k\in
R\backslash\left\{  i\right\}  }a_{kj}=\sum_{j\in C}a_{ij}\sum_{k\in
R\backslash\left\{  i\right\}  }a_{kj}=\sum_{j\in C}\sum_{k\in R}a_{ij}%
a_{kj}-\sum_{j\in C}a_{ij}\\
&  =\sum_{j\in C}\sum_{k\in R}a_{ij}a_{kj}-r_{i}.
\end{align*}
Substituting the value of $\left\Vert A\left[  U,W\right]  \right\Vert $ in
(\ref{in3}), we see that
\[
\sum_{j\in C}\sum_{k\in R}a_{ij}a_{kj}-\left(  \left(  s-k-1\right)
^{1/t}m^{1-1/t}+k\right)  r_{i}+\left(  t-1\right)  m^{1+\left(  k-1\right)
/t}\leq0.
\]
Summing this inequality for all $i\in R,$ we get%
\[
\sum_{i\in R}\sum_{j\in C}\sum_{k\in R}a_{ij}a_{kj}-\left(  \left(
s-k-1\right)  ^{1/t}m^{1-1/t}+k\right)  \left\Vert A\right\Vert +\left(
t-1\right)  m^{2+\left(  k-1\right)  /t}\leq0.
\]
Now note that%
\[
\sum_{i\in R}\sum_{j\in C}\sum_{k\in R}a_{ij}a_{kj}=\sum_{j\in C}\sum_{i\in
R}\sum_{k\in R}a_{ij}a_{kj}=\sum_{j\in C}r_{i}^{2}\geq\frac{1}{n}\left\Vert
A\right\Vert ^{2},
\]
and so,
\[
\frac{1}{n}\left\Vert A\right\Vert ^{2}-\left(  \left(  s-k-1\right)
^{1/t}m^{1-1/t}+k\right)  \left\Vert A\right\Vert -\left(  t-1\right)
m^{2+\left(  k-1\right)  /t}\leq0.
\]
Solving this inequality, we find that%
\[
\left\Vert A\right\Vert \leq\left(  1+\sqrt{1+\frac{4\left(  t-1\right)
m^{2+\left(  k-1\right)  /t}}{n\left(  \left(  s-k-1\right)  ^{1/t}%
m^{1-1/t}+k\right)  ^{2}}}\right)  \frac{\left(  \left(  s-k-1\right)
^{1/t}m^{1-1/t}+k\right)  n}{2}%
\]
and bounding the radical by the Bernoulli inequality, we obtain
\begin{align*}
\left\Vert A\right\Vert  &  \leq\left(  1+1+\frac{2\left(  t-1\right)
m^{2+\left(  k-1\right)  /t}}{n\left(  \left(  s-k-1\right)  ^{1/t}%
m^{1-1/t}+k\right)  ^{2}}\right)  \frac{\left(  \left(  s-k-1\right)
^{1/t}m^{1-1/t}+k\right)  n}{2}\\
&  =\left(  s-k-1\right)  ^{1/t}nm^{1-1/t}+kn+\frac{\left(  t-1\right)
m^{2+\left(  k-1\right)  /t}}{\left(  s-k-1\right)  ^{1/t}m^{1-1/t}+k}\\
&  \leq\left(  s-k-1\right)  ^{1/t}nm^{1-1/t}+\left(  t-1\right)
m^{1+k/t}+kn.
\end{align*}
This completes the induction step and the proof of Theorem \ref{thZ}.
\end{proof}

Stated in terms of bipartite graphs, Theorem \ref{thZ} is equivalent to the
following one:

\begin{theorem}
\label{thBG}Let $s\geq2,$ $t\geq2,$ $0\leq k\leq s-2,$ and let $G\left(
A,B\right)  $ be a bipartite graph with parts $A$ and $B.$ Suppose that $G$
contains no copy of $K_{s,t}$ with a vertex class of size $s$ in $A$ and a
vertex class of size $t$ in $B.$ Then $G\left(  A,B\right)  $ has at most
\[
\left(  s-k-1\right)  ^{1/t}\left\vert B\right\vert \left\vert A\right\vert
^{1-1/t}+\left(  t-1\right)  \left\vert A\right\vert ^{1+k/t}+k\left\vert
B\right\vert
\]
edges.
\end{theorem}

\subsection*{The proof of Theorem \ref{th1}}

\paragraph{Some graph notation}

Our graph notation follows \cite{Bol98}; in particular, given a graph $G$ and
a vertex $u$ of $G,$ we write:\smallskip

- $V\left(  G\right)  $ for the vertex set of $G;$

- $E\left(  G\right)  $ for the edge set of $G$ and $e\left(  G\right)  $ for
$\left\vert E\left(  G\right)  \right\vert ;$

- $G-u$ for the graph obtained from $G$ by removing the vertex $u.$

- $\Gamma\left(  u\right)  $ for the set of neighbors of $u$ and $d\left(
u\right)  $ for $\left\vert \Gamma\left(  u\right)  \right\vert .$\bigskip

\begin{proof}
[\textbf{Proof of Theorem \ref{th1}}]Inequality (\ref{in0}) has been proved in
\cite{Nik07d}, so we shall assume that $s\geq3$ and $t\geq3.$ Let $u\in
V\left(  G\right)  $ be any vertex of $G,$ let $U$ and $W$ be disjoint sets
satisfying $\left\vert U\right\vert =d\left(  v\right)  $ and $\left\vert
W\right\vert =n-1,$ and let $\varphi_{U}$ and $\varphi_{W}$ be bijections
\[
\varphi_{U}:U\rightarrow\Gamma\left(  u\right)  ,\text{ \ \ }\varphi
_{W}:W\rightarrow V\left(  G\right)  \backslash\left\{  u\right\}  .
\]
Define a bipartite graph $H$ with vertex classes $U$ and $W$ by joining $v\in
U$ and $w\in W$ whenever $\left\{  \varphi_{U}\left(  v\right)  ,\varphi
_{W}\left(  w\right)  \right\}  \in E\left(  G\right)  .$

We claim that $H$ does not contain a copy of $K_{t,s-1}$ with $s-1$ vertices
in $W$ and $t$ vertices in $U.$ Indeed, the map $\psi:V\left(  H\right)
\rightarrow V\left(  G\right)  $ defined as%
\[
\psi\left(  x\right)  =\left\{
\begin{array}
[c]{cc}%
\varphi_{U}\left(  x\right)  & \text{if }x\in U\\
\varphi_{W}\left(  x\right)  & \text{if }x\in W
\end{array}
\right.
\]
is a homomorphism of $H$ into $G-v.$ Assume for a contradiction that $F\subset
H$ is a copy of $K_{t,s-1}$ with a set $S$ of $s-1$ vertices in $W$ and a set
$T$ of $t$ vertices in $U$. Clearly $S$ and $T$ are the vertex classes of $F.$
Note that $\psi\left(  F\right)  $ is a copy of $K_{t,s-1}$ in $G-u,$ and
$\psi\left(  T\right)  =\varphi_{U}\left(  T\right)  \subset\Gamma_{G}\left(
u\right)  $ is the vertex class of $\psi\left(  F\right)  $ of size $t;\ $now,
adding $u$ to $\psi\left(  F\right)  ,$ we see that $G$ contains a $K_{t,s},$
a contradiction proving the claim.

Suppose that $0\leq k\leq\min\left\{  s,t\right\}  -2.$ Setting $k^{\prime
}=k-1,$ $s^{\prime}=s-1,$ $t^{\prime}=t,$ $A=W,$ $B=U,$ Theorem \ref{thBG}
implies that
\begin{align*}
e\left(  H\right)   &  \leq\left(  s-k-1\right)  ^{1/t}\left\vert U\right\vert
\left\vert W\right\vert ^{1-1/t}+\left(  k-1\right)  \left\vert U\right\vert
+\left(  t-1\right)  \left\vert W\right\vert ^{1+\left(  k-1\right)  /t}\\
&  \leq\left(  s-k-1\right)  ^{1/t}d\left(  u\right)  n^{1-1/t}+\left(
k-1\right)  d\left(  u\right)  +\left(  t-1\right)  n^{1+\left(  k-1\right)
/t}.
\end{align*}
On the other hand, we see that%
\[
e\left(  H\right)  =\sum_{v\in\Gamma\left(  u\right)  }d\left(  v\right)
-d\left(  u\right)  ,
\]
and so,
\begin{equation}
\sum_{v\in\Gamma\left(  u\right)  }d\left(  v\right)  \leq\left(  \left(
s-k-1\right)  ^{1/t}n^{1-1/t}+k\right)  d\left(  u\right)  +\left(
t-1\right)  n^{1+\left(  k-1\right)  /t}.\label{sq2}%
\end{equation}
Letting $A$ be the adjacency matrix of $G,$ note that the $u$'th row sum of
the matrix
\[
C=A^{2}-\left(  \left(  s-k-1\right)  ^{1/t}n^{1-1/t}+k\right)  A
\]
is equal to
\[
\sum_{v\in\Gamma\left(  u\right)  }d\left(  v\right)  -\left(  \left(
s-k-1\right)  ^{1/t}n^{1-1/t}+k\right)  d\left(  u\right)  ;
\]
consequently, the maximum row sum $r_{\max}$ of $C$ satisfies%
\[
r_{\max}\leq\left(  t-1\right)  n^{1+\left(  k-1\right)  /t}.
\]
Letting $\mathbf{x}$ be an eigenvector of $A$ to $\mu,$ we see that the value%
\[
\lambda=\mu^{2}-\left(  \left(  s-k-1\right)  ^{1/t}n^{1-1/t}+k\right)  \mu
\]
is an eigenvalue of $C$ with eigenvector $\mathbf{x}$. Therefore,%
\[
\mu^{2}-\left(  \left(  s-k-1\right)  ^{1/t}n^{1-1/t}+k\right)  \mu
=\lambda\leq r_{\max}\leq\left(  t-1\right)  n^{1+\left(  k-1\right)  /t}.
\]
Solving this inequality we obtain%
\begin{align*}
\mu &  \leq\left(  1+\sqrt{1+\frac{4\left(  t-1\right)  n^{1+\left(
k-1\right)  /t}}{\left(  \left(  s-k-1\right)  ^{1/t}n^{1-1/t}+k\right)  ^{2}%
}}\right)  \frac{\left(  s-k-1\right)  ^{1/t}n^{1-1/t}+k}{2}\\
&  \leq\left(  1+1+\frac{2\left(  t-1\right)  n^{1+\left(  k-1\right)  /t}%
}{\left(  \left(  s-k-1\right)  ^{1/t}n^{1-1/t}+k\right)  ^{2}}\right)
\frac{\left(  s-k-1\right)  ^{1/t}n^{1-1/t}+k}{2}\\
&  \leq\left(  s-k-1\right)  ^{1/t}n^{1-1/t}+\left(  t-1\right)  n^{k/t}+k.
\end{align*}
Now, if $s\geq t\geq3$, setting $k=t-2$, we obtain inequality (\ref{in1}),
completing the proof of Theorem \ref{th1}.
\end{proof}

\paragraph*{\textbf{Acknowledgement}}

Thanks are due to L\'{a}szl\'{o} Babai for details on his work with B.
Guiduli, and to Tibor Szab\`{o} for pointing out the relevance of norm-graphs
to the present topic.

\end{document}